\documentclass[a4paper,11pt,reqno]{amsart}

\usepackage{times}
\usepackage[T1]{fontenc}
\usepackage{mathrsfs}

\numberwithin{equation}{section}

\newtheorem{theorem}{Theorem}[section]

\newtheorem{corollary}[theorem]{Corollary}
\newtheorem{lemma}[theorem]{Lemma}
\newtheorem{proposition}[theorem]{Proposition}
\theoremstyle{remark}

\theoremstyle{definition}

\newcommand\bp{\begin{proof}}
\newcommand\ep{\end{proof}}

\newcommand\aaa{\mathfrak a}
\newcommand\bb{\mathfrak b}
\newcommand\mm{\mathfrak m}
\newcommand\pp{\mathfrak p}
\newcommand\qq{\mathfrak q}

\newcommand{\C}{\mathbb C}
\newcommand{\F}{\mathbb F}
\newcommand{\N}{\mathbb N}
\newcommand{\Q}{\mathbb Q}
\newcommand{\R}{\mathbb R}
\newcommand{\Z}{\mathbb Z}

\newcommand\A{\mathbb{A}}
\newcommand\ak{{\mathbb A}_K}
\newcommand\aks{{\mathbb A}_{S(\mm)}}
\newcommand\akf{{\mathbb A}_{K,f}}
\newcommand\oas{\hat{\mathcal O}_{S(\mm)}}
\newcommand\OO{{\mathcal O}}
\newcommand\ohs{{\hat{\OO}^*}}

\newcommand\Gal{\operatorname{Gal}}
\newcommand{\RR}{\mathcal R}
\newcommand{\PP}{\mathcal P}

\newcommand\eps{\varepsilon}

\newcommand{\Pm}{{\mathcal{P}_{\mathfrak{m}, 1}}}

\newcommand{\Kp}{K^{+}}
\newcommand{\bv}{{\bf v}}
\newcommand{\bw}{{\bf w}}
\newcommand{\bx}{{\bf x}}
\newcommand{\gchar}{{\chi}}

\begin{document}

\title{Ergodicity of the action of $K^*$ on $\ak$}

\author[J.C. Lagarias]{Jeffrey C. Lagarias}

\address{Department of Mathematics, University of Michigan, Ann Arbor, MI 48109-1043, USA}

\email{lagarias@umich.edu}

\author[S. Neshveyev]{Sergey Neshveyev}

\address{Department of Mathematics, University of Oslo,
P.O. Box 1053 Blindern, N-0316 Oslo, Norway}

\email{sergeyn@math.uio.no}

\thanks{Research supported in part by NSF grant DMS-1101373 and the Norwegian Research Council.}

\begin{abstract}
Connes gave a spectral interpretation of the critical zeros of zeta- and $L$-functions for a global field~$K$ using a space of square integrable
functions on the space $\ak/K^*$ of adele classes. It is known that for $K=\Q$ the space $\ak/K^*$ cannot be understood classically, or in other words, the action of $\Q^*$ on $\A_\Q$ is ergodic. We prove that the same is true for any global field $K$, in  both  the number field and function field cases.
\end{abstract}

\date{November 14, 2012; minor changes June 8, 2013: v. 2.7}

\maketitle

\section*{Introduction}

For a global field $K$, the space of adeles $\ak$ appeared explicitly (as ``valuation vectors'') in Tate's thesis~\cite{Ta50}. The action of $K^*$ on $\ak$ played a key role in his proof of the functional equation for Hecke $L$-functions. More recently, this action, together with the associated space $\ak/K^*$ of adele classes, was used by Connes~\cite{Con} to construct a Hilbert space that provides a spectral interpretation of the critical zeros of Hecke $L$-functions. One of the major difficulties in constructing such a Hilbert space, as well as in trying to prove the Riemann hypotheses following Connes' approach, is that the quotient space $\ak/K^*$ is quite tricky. The precise statement for $K=\Q$ is that the action of $\Q^*$ on $\A_\Q$ is ergodic, so from the point of view of classical measure theory the space $\A_\Q/\Q^*$ consists of just one point. This is a known byproduct of the analysis of the Bost-Connes system~\cite{bos-con,N1}. The aim of this paper is to prove that the same ergodicity holds for $K^*$ acting on~$\A_K$  for an arbitrary global field~$K$, i.e., for $K$ an algebraic number field or a one-variable algebraic function field over a finite field.

Partial results in this direction have been already known from the analysis of Bost-Connes type systems for global fields~\cite{LLN,N2,NR}. For example, for number fields, by \cite[Corollary~3.2]{N2}, the action of the subgroup of totally positive elements in $K^*$ on $\R^*_+\times\akf$, given by $x(t,a)=(N((x))t,xa)$, is ergodic. This result could be used to simplify some of our arguments. However, the simplifications would be nonessential,
and roughly amount to reducing the arguments below to the case of Hecke characters ($\bmod\ 1$)  instead of using all Hecke characters. We have therefore chosen to develop all the necessary tools from scratch. As a consequence, our arguments can be used to give different proofs of several results in~\cite{bos-con,LLN,N2,NR}, see Sections~\ref{s4.4} and~\ref{s5func}.

\bigskip


\section{Main result} \label{s1}

Let $K$ be a global field. Denote by $\ak$ the adele ring of $K$ and $\mu_{\A}$ an additive Haar measure on $\ak$.

\begin{theorem} \label{tmain1}
The action of $K^*$ on $(\ak,\mu_\A)$  by multiplication is ergodic; that is, for every $K^*$-invariant measurable subset of $\ak$ either the subset itself or its complement has measure zero.
\end{theorem}

In order to shorten the notation we will usually omit explicit mention of the  measures. In particular, by ergodicity of an action we will always mean ergodicity with respect to an appropriate Haar measure, or with respect to the push-forward of such a measure, unless said otherwise. The ring $\ak$ is not compact, and any choice of additive Haar measure $\mu_{\A}$ for $\ak$ will assign infinite measure to the full space $\ak$. The Haar measure is a product-measure, and we assume for convenience that for every non-Archimedean place $v$ of~$K$ the Haar measure on $K_v$ is normalized so that the measure of the maximal compact subring $\OO_v\subset K_v$ is equal to $1$.

\smallskip

Our first goal is to replace $\ak$ by a more manageable space and the action of $K^*$ by an action of a group of ideals on it.

\smallskip

If $K$ is a number field, denote by $\OO$ the ring of integers in $K$. In the case of function fields we fix a place~$\infty$ of~$K$, and then denote by $\OO$ the ring of elements in $K$ that have no poles away from $\infty$, so $\OO=\cap_{v\ne\infty}(K\cap\OO_v)$.

For a nonzero ideal $\mm$ in $\OO$, denote by $I_\mm$ the group of fractional ideals in $\OO$ relatively prime to $\mm$. Denote by~$V_f$ the set of finite places of~$K$ and by $S(\mm)\subset V_f$ the set of prime divisors of $\mm$. Put
$$
\aks:=\underset{v\in V_f\setminus S(\mm)}{{\prod}'}K_v,\ \ \oas:=\prod_{v\in V_f\setminus S(\mm)}\OO_v,
$$
where $\prod^{'}$ denotes restricted direct product with respect to the $\OO_v$.
Writing $\mm$ as $\prod_{v\in S(\mm)}\pp_v^{m_v}$, we put
$$
W_\mm :=\oas^*\times\prod_{v\in S(\mm)}(1+\hat\pp_v^{m_v})\subset\ohs:=\prod_{v\in V_f}\OO_v^*,
$$
where $\hat\pp_v$ is the maximal ideal in $\OO_v$. Identifying $I_\mm$ with $\aks^*/\oas^*$ we get actions of $I_\mm$ on $\ak^*/W_\mm$ and $\aks/\oas^*$ by multiplication. The Haar measure on $\aks$, normalized so that the measure of $\oas$ is equal to~$1$, defines a measure on $\aks/\oas^*$. The quotient $\ak^*/K^*W_\mm$ of the idele class group is a locally compact group under multiplication, so it also carries a Haar measure.

\begin{theorem} \label{tmain2}
The diagonal action of $I_\mm$ on $(\ak^*/K^*W_\mm)\times(\aks/\oas^*)$, equipped with the product measure, is ergodic for every nonzero ideal $\mm\subset\OO$.
\end{theorem}

In the remaining part of this section we will show that the statements of Theorems~\ref{tmain1} and~\ref{tmain2} are equivalent.
In the next two sections we  then prove Theorem \ref{tmain2} for number fields. In Section~\ref{s2} we give a preliminary result on asymptotic range of measurable cocycles, and in Section~\ref{s3} we establish Theorem \ref{tmain2}. In Section~\ref{sf} we describe  modifications to the proof needed to treat the function field case. Finally, in Section~\ref{s4}, we give various remarks.

\smallskip

The  compact groups $W_\mm$ form a base of neighbourhoods of the identity in $\ohs$. Hence, by averaging over~$W_\mm$, we can approximate any bounded measurable function on $\ak$ by a $W_\mm$-invariant function for some~$\mm$. It follows that the action of $K^*$ on $\ak$ is ergodic if and only if the action of $K^*$ on $\ak/W_\mm$ is ergodic for all $\mm$.

Denote by $K^*_{\mm,1}\subset K^*$ the subgroup of elements $a\equiv1\bmod\ \mm$, and by $\Pm \subset I_\mm$ the subgroup of principal fractional ideals
having a generator in  $K^*_{\mm,1}$. Put $U_{\mm,1}:=K^*_{\mm,1}\cap\OO^*$. Consider the subset
$$
X_\mm :=\left(K^*_\infty\times\aks\times\left(\prod_{v\in S(\mm)}(1+\hat\pp_v^{m_v})\right)\right)/W_\mm\subset\ak/W_\mm.
$$
Here, if $K$ is a number field, then $K_\infty=\prod_{v|\infty}K_v=\R^{r_1} \times \C^{r_2}$
(see Section \ref{s3}), and if $K$ is a function field, then $K_\infty$ is the completion of~$K$ corresponding to our fixed place $\infty$. The set $X_\mm$ is $K^*_{\mm,1}$-invariant, the set $K^*X_\mm\subset\ak/W_\mm$ is a subset of full measure, and $x X_\mm\cap X_\mm=\emptyset$ for any $x\in K^*\setminus K^*_{\mm,1}$. It follows that the action of $K^*$ on $\ak/W_\mm$ is ergodic if and only if the action of~$K^*_{\mm,1}$ on
$$
X_\mm\cong K^*_\infty\times\aks/\oas^*
$$
is ergodic, or in other words, the action of $\Pm\cong K^*_{\mm,1}/U_{\mm,1}$ on $K^*_\infty/U_{\mm,1}\times\aks/\oas^*$ is ergodic.

Now, for each $\mm$, induce the last action of $\Pm$ to an action of $I_\mm$. So consider the space
$$
Y_{\mm}:= I_\mm\times_{\Pm}(K^*_\infty/U_{\mm,1}\times\aks/\oas^*),
$$
with the action of $I_\mm$ defined  by its action on $I_{\mm}$ in the first factor by translations.
This action is ergodic if and only if the action of
$\Pm$ on $K^*_\infty/U_{\mm,1}\times \aks/\oas^*$ is ergodic. The map
$$
I_\mm\times K^*_\infty/U_{\mm,1}\times \aks/\oas^* \to I_\mm\times K^*_\infty/U_{\mm,1}\times\aks/\oas^*,\ \ (\aaa,x,b)\mapsto(\aaa,x,\aaa b),
$$
defines a measure class preserving homeomorphism of $Y_\mm$ onto
$$
Z_{\mm} :=(I_\mm\times_{\Pm} K^*_\infty/U_{\mm,1})\times \aks/\oas^*.
$$
Under this homeomorphism the action of $I_\mm$ on $Y_{\mm}$ becomes the diagonal action on $Z_{\mm}$, with the action of $I_\mm$ on $I_\mm\times_{\Pm} K^*_\infty/U_{\mm,1}$ defined  by the action on $I_{\mm}$ in the first factor by multiplication.

Finally, the group $I_\mm\times_{\Pm} K^*_\infty/U_{\mm,1}$ is isomorphic to the quotient of $\aks^*/\oas^*\times K^*_\infty$ by the subgroup $\{(x^{-1},x)\mid x\in K^*_{\mm,1}\}$. Applying the automorphism $(x,y)\mapsto(x,y^{-1})$ of $\aks^*/\oas^*\times K^*_\infty$ we conclude that $I_\mm\times_{\Pm} K^*_\infty/U_{\mm,1}$ is isomorphic to the quotient of $\aks^*/\oas^*\times K^*_\infty$ by $K^*_{\mm,1}$ under the diagonal embedding. The group $\aks^*/\oas^*\times K^*_\infty$ can be identified with
$$
H_\mm:=\left(K^*_\infty\times\aks^*\times\left(\prod_{v\in S(\mm)}(1+\hat\pp_v^{m_v})\right)\right)/W_\mm\subset\ak^*/W_\mm.
$$
We have $K^*H_\mm=\ak^*/W_\mm$ and $x H_\mm\cap H_\mm=\emptyset$ for any $x\in K^*\setminus K^*_{\mm,1}$. It follows that
$$
H_\mm/K^*_{\mm,1}\cong \ak^*/K^*W_\mm.
$$
To summarize, we have an isomorphism
$$
I_\mm\times_{\Pm} K^*_\infty/U_{\mm,1}\cong  \ak^*/K^*W_\mm
$$
that respects the embeddings of $I_\mm\cong \aks^*/\oas^*$. Therefore the action of $K^*$ on $\ak/W_\mm$ is ergodic if and only if the action of on $I_\mm$ on
$Z_\mm\cong (\ak^*/K^*W_\mm)\times(\aks/\oas^*)$ is ergodic.

This holds for all $\mm$, so the statements of Theorems~\ref{tmain1} and~\ref{tmain2} are indeed equivalent.

\bigskip

\section{Asymptotic range of a cocycle}\label{s2}

Let $(X,\mu)$ be a standard measure space, meaning that by removing a set of measure zero we can identify~$X$ with a complete separable metric space, in such a way that $\mu$ becomes the completion of a $\sigma$-finite Borel measure. Let $\RR$ be an ergodic countable equivalence relation on $(X,\mu)$, as in Feldman and Moore~\cite{FM1}. Recall that {\em ergodicity} means if a measurable subset of $X$ is a union of equivalence classes, then either the subset itself or its complement has measure zero. Assume $c$ is a $1$-cocycle on $\RR$ with values in an abelian second countable locally compact group $\Gamma$, so $c$ is a measurable map $\RR\to\Gamma$ such that $c(x,y)c(y,z)=c(x,z)$ for $x\sim_\RR y \sim_\RR z$. Associate to $c$ an equivalence relation $\RR(c)$ on $\Gamma\times X$ defined by
$$
(\beta,x)\sim_{\RR(c)}(\gamma,y)\ \ \text{iff}\ \ x\sim_\RR y\ \ \text{and}\ \ c(x,y)\beta=\gamma.
$$

For a measurable subset $A\subset X$ of positive measure, denote by $c_A$ the restriction of $c$ to $\RR\cap(A\times A)$. The {\em asymptotic range} $r^*(c)\subset\Gamma$ of $c$ is defined as the intersection of essential ranges of $c_A$ for all measurable subsets $A\subset X$ of positive measure~\cite[Definition~8.2]{FM1}. In other words, $\gamma\in \Gamma$ belongs to $r^*(c)$ if and only if for every neighbourhood $U$ of $\gamma$ and every measurable subset $A\subset X$ of positive measure there exists a measurable subset $B\subset A$ of positive measure and a measurable map $T\colon B\to A$ with graph in $\RR$ such that $c(x,Tx)\in U$ for all $x\in B$. We recall the following facts.

\begin{proposition} \label{asymprange}
The asymptotic range has the following properties:
\begin{enumerate}
\item[(i)]
the set $r^*(c)\subset \Gamma$ is a closed subgroup;
\item[(ii)]
for any measurable subset $A\subset X$ of positive measure we have $r^*(c_A)=r^*(c)$;
\item[(iii)]
the equivalence relation $\RR(c)$ on $\Gamma\times X$ is ergodic if and only if $r^*(c)=\Gamma$.
\end{enumerate}
\end{proposition}

\bp
Properties (i) and (ii) easily follow from definitions. Property (iii) is nontrivial, it follows from Theorem 8 in \cite{FM1} (see also Propositions~8.1 and~8.3 there).
\ep

Consider now a particular class of ergodic equivalence relations. Take a sequence $\{(X_n,\mu_n)\}^\infty_{n=1}$ of at most countable probability spaces. Put
$$
(X,\mu):=\prod^\infty_{n=1}(X_n,\mu_n).
$$
Define an equivalence relation $\RR$ on $X$ by
$$
x\sim_\RR y\ \ \text{iff}\ \ x_n=y_n\ \ \text{for almost all}\ \ n.
$$

Assume we are given a $\Gamma$-valued $1$-cocycle $c$ on $\RR$. An important example is the Radon-Nikodym cocycle~$c_\mu$. It is the $\R^*_+$-valued cocycle defined by
$$
c_\mu(x,y)=\prod^\infty_{n=1}\frac{\mu_n(y_n)}{\mu_n(x_n)}.
$$
More generally, we say that a $\Gamma$-valued cocycle $c$ is of {\em product type} if there exist maps $\theta_n\colon X_n\to\Gamma$ such that
$$
c(x,y)=\prod^\infty_{n=1}\theta_n(y_n)\theta_n(x_n)^{-1}.
$$
It is not difficult to see that a cocycle $c$ is of product type if and only if it is constant on the sets $$\{(x,y)\mid x_n=a,\ y_n=b\ \text{and}\ x_m=y_m\ \text{for}\ m\ne n\}\subset\RR$$ for all $n\in\N$ and $a,b\in X_n$.

For a finite subset $I\subset\N$ and $a\in\prod_{n\in I}X_n$ put
$$
Z(a):=\{x\in X\mid x_n=a_n\ \text{for}\ n\in I\}.
$$
Motivated by~\cite{AW}, define the {\em asymptotic ratio set} $r(c)$ of a $\Gamma$-valued cocycle $c$ as the set of elements $\gamma\in\Gamma$ such that for every neighbourhood $U$ of $\gamma$ there exist numbers $t>s>0$, a sequence $\{I_n\}^\infty_{n=1}$ of mutually disjoint finite subsets of~$\N$, subsets
$K_n,L_n\subset\prod_{k\in I_n}X_k$ and bijections $T_n\colon
K_n\to L_n$ such that
$$
c_\mu(x,\bar T_nx)\in[s,t]\ \text{and}\ c(x,\bar T_nx)\in U\ \text{for all}\
x\in Z(a),\ a\in K_n\ \text{and}\ n\ge1,
$$
where $\bar T_n$ denotes the obvious map $\cup_{a\in K_n}Z(a)\to\cup_{b\in L_n}Z(b)$ defined by $T_n\colon K_n\to L_n$, and
$$
\sum_{n=1}^\infty\sum_{a\in K_n}\mu(Z(a))=\infty.
$$

\begin{proposition} \label{pratio}
For any $\Gamma$-valued $1$-cocycle $c$ on $\RR$ we have $r(c)\subset r^*(c)$. If the cocycle $c$ is of product type, and for every $\gamma\in\Gamma$ there exist  numbers $t>s>0$ and a neighbourhood $U$ of $\gamma$ such that $c_\mu(x,y)\in[s,t]$ as long as $c(x,y)\in U$, then $r(c)=r^*(c)$.
\end{proposition}

Note that for the existence of $s$ and $t$ in the above proposition it suffices to have a continuous homomorphism $p\colon \Gamma\to\R^*_+$ such that $p\circ c=c_\mu$. Similarly, if such $p$ exists, then the requirements on $s$, $t$ and $c_\mu$ in the definition of the asymptotic ratio set can be omitted.

\smallskip

Proposition~\ref{pratio} is well-known for the Radon-Nikodym cocycle $c_\mu$ and is considered to be obvious in this case~\cite{kr}. A rather detailed proof is contained in \cite[Appendix A]{NR}. Essentially the same arguments work in the general case. Let us only sketch a proof of the inclusion $r(c)\subset r^*(c)$, which is the one we are going to use later.

Assume $\gamma\in r(c)$. Fix a neighbourhood $U$ of $\gamma$. Choose $s$, $t$, $I_n$, $K_n$, $L_n$ and $T_n$ as in the definition of the asymptotic ratio set. Consider the sets $A_n:=\cup_{a\in K_n}Z(a)$ and $B_n:=\cup_{b\in L_n}Z(b)$, and put
$$
A:=\cup_n\left(A_n\setminus\cup_{m<n}(A_m\cup B_m)\right)\ \ \text{and}\ \ B:=\cup_n\left(B_n\setminus\cup_{m<n}(A_m\cup B_m)\right).
$$
Since the sets $A_n\cup B_n$ are mutually independent and
$$
\sum_n\mu(A_n\cup B_n)\ge\sum_n\mu(A_n)=\infty,
$$
we have $\mu(A\cup B)=1$. Define a measurable bijective map $T\colon A\to B$ as follows: if $$x\in (A_n\setminus\cup_{m<n}(A_m\cup B_m))\cap Z(a)$$ for some $a\in K_n$, then put $Tx:=\bar T_nx$. Then for any $x\in A$ we have
$c_\mu(x,Tx)\in[s,t]$ and $c(x, Tx)\in U$.

Observe now that if we fixed a cylindrical set $Z(a)$ and then used only $I_n$, $K_n$, $L_n$ and $T_n$ with sufficiently large~$n$ in the above construction of $A$, $B$ and $T$, then $T$ would map $Z(a)\cap A$ onto $Z(a)\cap B$. We can therefore conclude that if $Z\subset X$ is a finite union of cylindrical sets $Z(a)$, then there exist measurable subsets $A,B\subset Z$ and a measurable bijective map $T\colon A\to B$ with graph in~$\RR$ such that
$$
\mu(Z\setminus(A\cup B))=0,\ \ c_\mu(x,Tx)\in[s,t]\ \ \text{and}\ \ c(x,Tx)\in U\ \ \text{for}\ \ x\in A.
$$
Any measurable subset $Z'\subset X$ of positive measure can be approximated arbitrarily well by such~$Z$. Using that
$$
\frac{dT^{-1}\mu}{d\mu}(x)=c_\mu(x,Tx)\in[s,t]\ \ \text{for all}\ \ x\in A,
$$
it is easy to see that if $\mu(Z'\Delta Z)$ is sufficiently small, then $$\mu(Z'\cap T^{-1}(Z'\cap B))>0.$$
Therefore $T$ maps the subset $Z'\cap T^{-1}(Z'\cap B)$ of $Z'$ of positive measure into $Z'$, and $c(x,Tx)\in U$ for any~$x$ in this subset.
Since $U$ was arbitrary, this shows that $\gamma$ lies in the essential range of $c_{Z'}$.

%
%
\bigskip

\section{Distribution of prime ideals in number fields} \label{s3}

In this section we will prove Theorem~\ref{tmain2} for number fields $K$. Set
$n = [K:\Q] = r_1+2r_2$, so $\ak$ and~$\ak^{\ast}$ each have $r_1+ r_2$ Archimedean places, with $r_2$ of them complex. Write $x = (x_{v})_{v}$ for elements of~$\ak$, where $v$ runs over the valuations of $K$. Valuations are normalized as in Tate's thesis \cite[Sec. 2.1]{Ta50}  and denoted $||\cdot||_v$ following  Lang \cite[XIV, \S1]{La},  with normalized complex valuations defined via the complex norm $||z||_{\C} := |z|^2$. The idelic norm map $N \colon \ak^{*} \to \R^*_+$ sends $(x_v)_{v} \mapsto \prod_{v} ||x||_v$. For any nonzero integral ideal $\mm$, the  {\em quotient norm map}  $N_{\mm}: \ak^*/K^*W_\mm \to \R^{*}_+$ is well-defined, and we put
\[
\Gamma_{\mm} := \ker N_{\mm}  \subset \ak^*/K^*W_\mm .
\]
The norm map $N\colon \ak^*\to \R^*_+$ has a
(noncanonical)  splitting homomorphism $s\colon \R^*_+\to\ak^*$. We fix such a continuous homomorphism, for example, letting $s(t)_v=1$ for $v\in V_f$ and $s(t)= (t^{1/n}, t^{1/n}, \ldots , t^{1/n})$ at  the $r_1+r_2$ Archimedean places $v|\infty$ (note again that $||t^{1/n}||_{v} = t^{2/n}$ at complex places). Using this homomorphism we can identify $\ak^*/K^*W_\mm$ with $\R^*_+\times\Gamma_\mm$. Then the embedding
\[
I_\mm\hookrightarrow\R^*_+\times\Gamma_{\mm}=\ak^*/K^*W_\mm
\]
has the form $\aaa\mapsto(N(\aaa)^{-1},\rho_\mm(\aaa))$ for some homomorphism $\rho_\mm\colon I_\mm\to \Gamma_{\mm}$. Here $N(\aaa)$ denotes the norm of the ideal $\aaa$; note that if $x\in\akf^*$ is an idele representing $\aaa$, then $N(\aaa)=N(x)^{-1}$. The  group $\Gamma_{\mm}$ is a compact abelian Lie group, whose connected components have real dimension $n-1.$ The points $\rho_\mm(\aaa) \in \Gamma_{\mm}$ can be thought of as measuring {\em generalized angles} of  ideals $\aaa\in I_\mm$; see Sections \ref{s4.2} and \ref{s4.3} below.

Consider now the orbit equivalence relation defined by the action of $I_\mm$ on $\aks/\oas^*$. This action is essentially free, so outside a set of measure zero we can define an $\R^*_+\times\Gamma_{\mm}$-valued $1$-cocycle by
$$
c(x,y)=(N(\aaa),\rho_\mm(\aaa))\ \ \text{if}\ \ \aaa x=y.
$$
Furthermore, it is easy to see that the action of $I_\mm$ on $\aks/\oas^*$ is ergodic. (Further justification  is supplied in the proof of Theorem~\ref{tmain2} below; more precisely, we will see there that the orbit equivalence relation on $\oas/\oas^*$ is ergodic, and since every $I_\mm$-orbit in $\aks/\oas^*$ intersects $\oas/\oas^*$, it follows that the action of $I_\mm$ on $\aks/\oas^*$ is also ergodic.) Therefore the statement of Theorem~\ref{tmain2} is equivalent to the equality $r^*(c)=\R^*_+\times\Gamma_{\mm}$.

The computation of the asymptotic range of $c$ will rely on the following distribution result. A closely related result is \cite[Theorem XV.5.6]{La}, see also Section~5.2 below.
\begin{theorem}\label{tdistrib}
The image under $\rho_\mm$ of the set of prime ideals in $\OO$, ordered by the norm, is equidistributed in~$\Gamma_\mm$. More precisely, for any continuous function $f$ on $\Gamma_\mm$ we have
$$
\frac{1}{|\{\pp\notin S(\mm) : N(\pp)\le x\}|}\sum_{\pp\notin S(\mm):N(\pp)\le x}f(\rho_\mm(\pp))\to\int_{\Gamma_\mm}f\,d\lambda_\mm\ \ \text{as}\ \ x\to\infty,
$$
where $\lambda_\mm$ is the normalized Haar measure on the compact group $\Gamma_\mm$.
\end{theorem}

\bp In order to prove the theorem it suffices to show that for any nontrivial character $\chi$ of $\Gamma_\mm$ we have
$$
\frac{1}{|\{\pp\notin S(\mm) : N(\pp)\le x\}|}\sum_{\pp\notin S(\mm):N(\pp)\le x}\chi(\rho_\mm(\pp))\to0\ \ \text{as}\ \ x\to\infty.
$$
In other words, since $|\{\pp : N(\pp)\le x\}|\sim{x}/{\log x}$ by the prime ideal theorem, we have to show
\begin{equation} \label{emain}
\sum_{\pp\notin S(\mm) : N(\pp)\le x}(\chi\circ\rho_\mm)(\pp)=o\left(\frac{x}{\log x}\right).
\end{equation}
But by construction of $\rho_\mm$, the character $\chi\circ\rho_\mm$ of $I_\mm$ is a Hecke character ($\bmod\ \mm$) that is not of the form $\aaa\mapsto N(\aaa)^{it}$, $t\in\R$. For such characters \eqref{emain} is well-known, see e.g.\ Lang \cite[Theorem XV.5.5]{La} or Narkiewicz\ \cite[Prop.~7.17 (3rd Ed.);  Prop. 7.9 (2nd Ed.)]{Nark}.
\ep

\begin{corollary} \label{c1}
For any subset $A\subset\Gamma_\mm$ with boundary of measure zero, we have
$$
|\{\pp\notin S(\mm): N(\pp)\le x\ \text{and}\ \rho_\mm(\pp)\in A\}|=\lambda_\mm(A)\frac{x}{\log x}+o\left(\frac{x}{\log x}\right)\ \ \text{as}\ \ x\to\infty.
$$
\end{corollary}

\begin{corollary} \label{cmain}
For any subset $A\subset\Gamma_\mm$ with boundary of measure zero, and any $\delta>0$, we have
$$
|\{\pp\notin S(\mm): x<N(\pp)\le (1+\delta)x\ \text{and}\ \rho_\mm(\pp)\in A\}|=\lambda_\mm(A)\delta  \frac{x}{\log x}+o\left(\frac{x}{\log x}\right)\ \ \text{as}\ \ x\to\infty.
$$
\end{corollary}

\bp[Proof of Theorem~\ref{tmain2} for number fields] We are given an ideal $\mm$. Consider the subset $\oas/\oas^*\subset\aks/\oas^*$, which is  of measure one. Modulo a set of measure zero it can be identified with the product of spaces $\OO^\times_\pp/\OO_\pp^*\cong\Z_+$, for $\pp\notin S(\mm)$, where $\OO^\times_\pp=\OO_\pp\setminus\{0\}$ and $\Z_{+} = \Z_{\ge0}$ is the additive semigroup. Then the measure is the product of measures $\mu_\pp$ on $\Z_+$ defined by
$$
\mu_\pp(n):=N(\pp)^{-n}(1-N(\pp)^{-1}).
$$
In this picture the orbit equivalence relation on $\oas/\oas^*$ is exactly the ergodic equivalence relation considered in the previous section:
$$
a\sim b\ \ \text{if and only if} \ \ a_\pp=b_\pp\ \ \text{for almost all}\ \ \pp,
$$
and the cocycle $c_{\oas/\oas^*}$ is the cocycle of product type defined by the maps $\theta_\pp\colon \Z_+\to\R^*_+\times\Gamma_\mm$,
$$
\theta_\pp(n):=(N(\pp)^n,\rho_\mm(\pp)^n).
$$
Note that these formulas show in particular that if we define a homomorphism $p\colon\R^*_+\times\Gamma_\mm\to\R^*_+$ by $p(x,\gamma):=x^{-1}$, then $p\circ c$ coincides with the Radon-Nikodym cocycle.

In order to prove the theorem, by Proposition \ref{asymprange} it suffices to show that the asymptotic range of $c_{\oas/\oas^*}$ is equal to $\R^*_+\times\Gamma_\mm$. Fix a point $(x_0,y_0)\in \R^*_+\times\Gamma_\mm$. Since the asymptotic range is a closed subgroup, it suffices to consider $x_0>1$. By Proposition~\ref{pratio} it is enough to show that $(x_0,y_0)$ belongs to the asymptotic ratio set of~$c_{\oas/\oas^*}$. This is done along familiar lines~\cite{B,N2}, as follows.

Fix a neighbourhood $U$ of $(x_0,y_0)$. Choose $\eps>0$ and an open subset $V\subset \Gamma_\mm$ with boundary of measure zero such that
$$
(x_0-\eps,x_0+\eps)\times y_0VV^{-1}\subset U.
$$
Choose $\delta>0$ such that $1+\delta<x_0$ and $\delta x_0<\eps$. Define sets $B_n$ of prime ideals by
$$
B_{2k}=\{\pp\notin S(\mm): x_0^{2k}<N(\pp)\le(1+\delta)x_0^{2k}\ \text{and}\ \rho_\mm(\pp)\in V\},
$$
$$
B_{2k+1}=\{\pp\notin S(\mm): x_0^{2k+1}<N(\pp)\le(1+\delta)x_0^{2k+1}\ \text{and}\ \rho_\mm(\pp)\in y_0V\}.
$$
By the choice of $\delta$ the sets $B_n$ are pairwise disjoint. By Corollary~\ref{cmain} we have
$$
|B_n|\sim \frac{\lambda_\mm(V)\delta x^n_0}{n\log x_0}\ \ \text{as}\ \ n\to\infty.
$$
In particular, there exists $k_0$ such that $|B_{2k+1}|\ge |B_{2k}|$ for $k\ge k_0$. For every $k\ge k_0$ choose a subset $C_{2k+1}\subset B_{2k+1}$ such that $|C_{2k+1}|=|B_{2k}|$. Let $\pp_1,\pp_2,\dots$ be the elements of $\cup_{k\ge k_0}B_{2k}$ enumerated so that $N(\pp_1)\le N(\pp_2)\le\dots$, and let $\qq_1,\qq_2,\dots$ be the elements of $\cup_{k\ge k_0}C_{2k+1}$ enumerated in the same way. Then, if $\pp_n\in B_{2k}$ for some $k$, we have $\qq_n\in B_{2k+1}$, and therefore
$$
N(\qq_n\pp_n^{-1})\in(x_0-\eps,x_0+\eps)\ \ \text{and}\ \ \rho_\mm(\qq_n\pp_n^{-1})\in y_0VV^{-1}.
$$
We also have
$$
\sum_{n=1}^\infty N(\pp_n)^{-1}\ge\sum^\infty_{k=k_0}\frac{|B_{2k}|}{(1+\delta)x_0^{2k}}=\infty.
$$
Hence for the sets $I_n$, $K_n$ and $L_n$ required by the definition of the asymptotic ratio set we can take $I_n=\{\pp_n,\qq_n\}$, $K_n=\{(1,0)\}\subset\Z_+^{I_n}$ and $L_n=\{(0,1)\}$.
\ep

\bigskip

\section{Function fields}\label{sf}

In this section we will prove Theorem~\ref{tmain2} for function fields. The strategy is the same as for number fields, but some changes are necessary, since for function fields the Chebotarev density theorem does not hold for the natural density (the density oscillates as a function of $x$ with no limiting value), and therefore the straightforward analogue of Theorem~\ref{tdistrib} cannot be true.

\smallskip

Assume $K$ is a global function field with constant field $\F_q$. In this case the image of the norm map $N\colon\ak^*\to\R^*_+$ is $q^\Z$. Fix an idele $a$ of norm $q$. Then for a given nonzero ideal $\mm\subset\OO$ the group $\ak^*/K^*W_\mm$ can be identified with $a^\Z\times\Gamma_\mm\cong \Z\times\Gamma_\mm$, where as before $\Gamma_\mm$ is the kernel of $N_\mm\colon \ak^*/K^*W_\mm\to q^\Z$. We have to show that the diagonal action of~$I_\mm$ on $\Z\times\Gamma_\mm\times\aks/\oas^*$ is ergodic. Since the group $\OO^*_\infty$ is profinite, we can make one more reduction and replace $\Gamma_\mm$ by the quotient by a compact open subgroup $U\subset\OO^*_\infty$. The group $\Gamma_\mm/U$ is finite, and it is canonically isomorphic to the quotient of $\ak^*$ by the open subgroup $a^\Z K^*UW_\mm$. By class field theory it follows that $\Gamma_\mm/U$ is the Galois group of a finite abelian extension $L/K$.

\begin{lemma}
The extension $L/K$ is geometric; that is,  $K$ and $L$ have the same field
of constants.  It is unramified at all finite places $\pp\notin S(\mm)$.
\end{lemma}

\bp Assume the constant field of $L$ is $\F_{q^m}$. The Artin map $
\ak^*\to\Gal(L/K)$ defines a homomorphism $\ak^*\to\Gal(\F_{q^m}/\F_q)$.
This homomorphism factors through the divisor group of $K$ and maps a divisor~$D$ to the Frobenius raised to the power $\deg D$. But by definition of~$L$ our fixed idele $a$ acts trivially on $L$, so the corresponding divisor of degree $1$ acts trivially on $\F_{q^m}$. Hence $m=1$, so the extension is geometric. Since~$\oas^*$ is contained in the kernel of the homomorphism $\ak^*\to\Gal(L/K)$, the extension is unramified at all finite places $\pp\notin S(\mm)$.
\ep

Therefore in order to prove Theorem~\ref{tmain2} for function fields it suffices to establish the following.

\begin{theorem} \label{tmain3}
For any nonzero ideal $\mm\subset\OO$ and any finite abelian geometric extension $L/K$ unramified at all finite places $\pp\notin S(\mm)$, the action of $I_\mm$ on $\Z\times\Gal(L/K)\times\aks/\oas^*$, given by
\[
\aaa(n,g,x)=(n+\deg\aaa,\sigma(\aaa)g,\aaa x),
\]
where $\sigma\colon I_\mm\to\Gal(L/K)$ is the Artin map, is ergodic.
\end{theorem}

This result is, in fact, a consequence of \cite[Theorem~2.1]{NR}, but here we will sketch a slightly different proof
following the arguments of the previous section. It should also be remarked that the theorem is not true for
nongeometric extensions, see Section~\ref{s5func}.

\bp[Proof of Theorem~\ref{tmain3}] Similarly to the proof of Theorem~\ref{tmain2} for number fields, define a $\Z\times\Gal(L/K)$-valued cocycle $c$ on the orbit equivalence relation on $\aks/\oas^*$ by
$$
c(x,y)=(\deg\aaa,\sigma(\aaa))\ \ \text{if}\ \ \aaa x=y.
$$
We aim to show that the asymptotic range of $c$ is $\Z \times \Gal(L/K)$.
For this it suffices to show that for any $g\in\Gal(L/K)$ the element $(1,g)$ belongs to the asymptotic ratio set of the restriction of $c$ to the equivalence relation on $\prod_{\pp\notin S(\mm)}\OO^\times_\pp/\OO^*_\pp\cong\prod_{\pp\notin S(\mm)}\Z_+.$
For this, define sets $B_n$ of prime ideals by
$$
B_{2k}=\{\pp\notin S(\mm): \deg\pp=2k\ \text{and}\ \sigma(\pp)=e\},
$$
$$
B_{2k+1}=\{\pp\notin S(\mm): \deg\pp=2k+1\ \text{and}\ \sigma(\pp)=g\}.
$$
By the Chebotarev density theorem for function fields, see e.g.\ \cite[Proposition~6.4.8]{FJ}, we have
$$
|B_n|=\frac{q^n}{n[L:K]}+O(q^{n/2})\ \ \text{as}\ \ n\to\infty.
$$
Proceeding now as in the proof of Theorem~\ref{tmain2} for number fields, we conclude that $(1,g)$ belongs to the asymptotic ratio set.
\ep


\bigskip

\section{Remarks}\label{s4}

\subsection{} Ergodicity of an action implies that the orbit of almost every point is dense. However, the  converse assertion  is not true in general.
Here it is not difficult to  show directly that  almost every orbit
of the action of~$K^*$ on $\ak$ is dense. More precisely, we show that the orbit of every noninvertible adele $a$ such that $a_v\ne0$ for all places~$v$ is dense. This was observed for $K=\Q$ by Laca and Raeburn \cite[Lemma~3.2]{LR}. In order to prove the same for an arbitrary global field $K$, fix $b\in\ak$ and a neighbourhood~$U$ of~$b$. We want to find $x\in K^*$ such that $xa\in U$. The set $U$ contains an open subset of the form $V\times (b_f+\hat\aaa)$, where $V\subset K_\infty\setminus\{0\}$, $\aaa$ is a nonzero ideal in~$\OO$ and $\hat\aaa$ is the closure of $\aaa$ in $\hat\OO$.

We claim that there exists $q\in K^*$ such that $qa_f\in\hat\aaa$ and $a_\infty^{-1}V+q\OO=K_\infty$. In order to see this, choose $y\in\OO^\times$ such that $ya_f\in\hat\aaa$. Consider the set $\{v_1,v_2,\dots\}$ of finite places $v\notin S(\aaa)$ such that $\|a_v\|_v<1$. Note that since $a$ is not invertible and $a_v\ne0$ for all $v$, this set is infinite. Consider the sequence of ideals $\mm_n=\pp_{v_1}\dots\pp_{v_n}$ in $\OO$. Choose a subsequence $\{\mm_{n_k}\}^\infty_{k=1}$ belonging to the same ideal class. Then the ideals $\mm_{n_1}^{-1}\mm_{n_k}=\pp_{v_{n_1+1}}\pp_{v_{n_1+2}}\dots\pp_{v_{n_k}}$ are principal. Let $z_k$ be a generator of $\mm_{n_1}^{-1}\mm_{n_k}$. Then the elements $q_k=z_k^{-1}y\in K^*$ have the properties that $q_ka_f\in\hat\aaa$ and $N((q_k))\to0$ as $k\to\infty$. Since~$\OO$ is a cocompact lattice in $K_\infty$ and since for appropriate units $u_k\in\OO^*$ we have $\|u_kq_k\|_v\to0$ for every~$v|\infty$, for sufficiently large $k$ we get $a_\infty^{-1}V+q_k\OO=K_\infty$. Thus our claim is proved.

We can now find the required element $x$. First, using density of $K$ in $\akf$, choose $x'\in K$ such that $x'a_f\in b_f+\hat\aaa$. Next, using that $a_\infty^{-1}V+q\OO=K_\infty$, choose $y\in\OO$ such that $x'+qy\in a_\infty^{-1}V$. Then the element $x=x'+qy\in K$ has the properties that $xa_\infty\in V$ and $xa_f=x'a_f+qya_f\in b_f+\hat\aaa$, as $qa_f\in\hat\aaa$. Therefore $xa\in V\times (b_f+\hat\aaa)\subset U$. Since $0\notin V$, we must have $x\in K^*$.

\subsection{}\label{s4.2}
For general number fields $K$  the key result underlying the proof of Theorem~\ref{tmain1} is Theorem~\ref{tdistrib}, which gives  equidistribution over the group $\Gamma_\mm$. This result  uses the full set of Hecke characters (``Gr\"{o}ssencharacters'') of the field~$K$, as introduced in Hecke  \cite{he18}, \cite{he20}, and given in Narkiewicz \cite[Sect. 7.3 (2nd Ed)]{Nark} or Neukirch \cite[Chap. VII.6]{Neu}. If $K \ne \Q$, then the  set of such Gr\"{o}ssencharacters includes characters of infinite order. The equidistribution result in Theorem \ref{tdistrib} for number fields  is at bottom derived using  the fact that the $L$-functions of Gr\"{o}ssencharacters are free from zeros and poles on the line $\operatorname{Re}(s)=1$, except for the characters $N(\cdot)^{it}$ having simple poles at $s=1-it$. Since $L(s,N(\cdot)^{it}\chi)=L(s+it,\chi)$, it makes sense to normalize characters, that is, to fix representatives among the characters that differ by a factor $N(\cdot)^{it}$. In this and the next subsection we fix the section $s(t)=(t^{1/n},\dots,t^{1/n},1,1,\dots)$ of the norm map, and call a character of $\ak^*/K^*$ normalized if it is trivial on $s(\R^*_+)$. This normalization agrees with Lang~\cite{La}, but differs from Narkiewicz~\cite{Nark}.

In order to describe the nature of the  groups $\ak/K^*W_\mm$ and $\Gamma_{\mm}$ in Section \ref{s3} in
more classical language, observe that the quotient of $\ak/K^*W_\mm$ by $K^*_\infty$ is isomorphic to the finite group $C_\mm=I_\mm/\Pm$, the ray class group modulo $\mm$. We therefore have a short exact sequence
$$
1\to K_\infty^*/U_{\mm,1}\to \ak/K^*W_\mm\to C_\mm\to1.
$$
Hence $\ak/K^*W_\mm$ is an abelian Lie group with finitely many connected components, and with connected component of the identity isomorphic to $\Kp_\infty/U^+_{\mm,1}$, where $\Kp_\infty$ is the connected component of the identity in~$K^*_\infty$ and $U^+_{\mm,1}=U_{\mm,1}\cap\Kp_\infty$ is the group of totally positive units $\bmod\ \mm$.

For an Archimedean place $v$, put $N_v=1$ if $v$ is real and $N_v=2$ if $v$ is complex. Denote by $K^*_{\infty,1}\subset K^*_\infty$ the subgroup of elements $x$ such that $\prod_{v|\infty}||x_v||_v=\prod_{v|\infty}|x_v|^{N_v}=1$. Then we also have a short exact sequence
$$
1\to K^*_{\infty,1}/U_{\mm,1}\to \Gamma_\mm\to C_\mm\to1.
$$
The homomorphism $\rho_\mm\colon I_\mm\to\Gamma_\mm$ is such that the image of $\rho_\mm(\aaa)$ in $C_\mm$ is the class $[\aaa]$ of $\aaa$, and if $\aaa=(x)$ for some $x\in K^*_{\mm,1}$, then
$$
\rho_\mm(\aaa)=s_\infty(N(\aaa))x^{-1}U_{\mm,1}\in K^*_{\infty,1}/U_{\mm,1}\subset\Gamma_\mm,
$$
where $s_\infty(t)=(t^{1/n},\dots,t^{1/n})\in K^*_\infty$. More explicitly, the group $\Gamma_\mm$ can be described as the quotient of $I_\mm\times K^*_{\infty,1}/U_{\mm,1}$ by $\Pm$ under the diagonal embedding, where the homomorphism $\Pm\to K^*_{\infty,1}/U_{\mm,1}$ is given by $(x)\mapsto s_\infty(N((x)))^{-1}xU_{\mm,1}$ for  $x\in K^*_{\mm,1}$.

Let $r_1$ and $r_2$ be the number of real and complex places of $K$, respectively, with $n=r_1+2r_2$.  Then $K_{\infty} = \R^{r_1} \times \C^{r_2}$. The connected component $\Kp_\infty/U^+_{\mm,1}$ of the identity in $\ak/K^*W_\mm$ has the form $\R^{r_1+2r_2}/\Lambda_\mm$, for a certain lattice $\Lambda_\mm$ of (nonfull) rank $r_1+2r_2- 1$, which  (under a logarithmic change of variable) encodes  the unit lattice together with additional lattice vectors encoding the ambiguity of the argument at each complex place. The connected component of the identity of the smaller group $\Gamma_{\mm}$ is then $V/\Lambda_\mm$, where $V$ is the real subspace of $\R^n$ of dimension $n-1$ spanned by the lattice~$\Lambda_\mm$. We emphasize that the lattice $\Lambda_m$ has  rank $n-1$, which for $r_2>0$ is larger than  the standard logarithmic  encoding in~$\R^{r_1+r_2}$ of the unit lattice of rank $r_1+r_2- 1$ used for computing the regulator of the field $K$ (as in Lang \cite[Chap. V]{La}).
The standard encoding detects only size information on the units, while the lattice $\Lambda_\mm$ has an extra $r_2$ components on which it encodes  information on the complex arguments (``angles'')  of the units.

We give an example. Take $K=\Q(\theta) $ to be a (non-Galois) real cubic field of discriminant $d_K= -23$ obtained by adjoining the real root  of $X^3-X-1=0$. This equation has one real root $\theta \approx 1.3247$ and two conjugate complex roots $\vartheta^{\pm} = \frac{1}{\sqrt{\theta}} e^{\pm2 \pi i  \phi}$. An important feature  is that the angle $\phi$ is irrational, because all nonzero powers of~$\vartheta^{+}$ are nonreal, since  $\Q(\vartheta^{+})$ contains no real units except $\pm1$. It is known that the field $K$ has  ring of integers $\OO = \Z[ 1, \theta, \theta^2].$ The unit rank of $K$ is $r_1+ r_2-1= 1$, and $\theta$ is a fundamental unit of $K$.  We choose $\mm= (1)$ to be the unit ideal. The field $K$ is  known to have class number $1$, so that $I_{\mm} = \Pm$, and  the group $\ak^*/K^*W_\mm=\ak^*/K^*\ohs$ is isomorphic to~$K^{*}_\infty/U_{\mm, 1}$, with $U_{\mm, 1} = U_K$ being the full unit group. Since $-1\in U_{\mm, 1}$, this group is connected, so it coincides with $\Kp_\infty/U^+_K$, where $U_K^+$ is the group of totally positive units, and $\Kp_{\infty} = \R_{+}^{*} \times \C^{*}$. Under the logarithmic change of coordinate sending $\log\colon \R^{*}_+ \to \R$ and
$\log\colon \C^{*} \to \R \oplus \R/2\pi \Z$, the group $U_K^+$  is mapped
to a $1$-dimensional lattice in $\R \oplus \R \oplus \R/2\pi \Z$, with basis vector:
$$
\bv_1= (\log \theta, - \frac{1}{2}\log \theta, 2\pi  \phi),
$$
with the argument defined $(\bmod~2\pi )$. We can lift  this lattice to $\R^3$, to obtain a $2$-dimensional lattice with  a second basis vector coming from the ambiguity of the argument:
$$
\bv_2=(0,0,  2 \pi ).
$$
Thus we obtain a $2$-dimensional lattice $\Lambda_{K}:=\Z [\bv_1, \bv_2]$ inside $\R^3$. Then the group $\ak^*/K^*\ohs=\Kp_\infty/U_K^+$ can be identified under the logarithmic map with $\R^3/ \Lambda_K$, with $\bv_2$ removing the ambiguity of the branch of the logarithm. Now let $V$ be the  $2$-dimensional subspace $\R[\bv_1, \bv_2]$, given explicitly by
$V=\{ \bx= (x_1, x_2, x_3)\, :\, x_1 + 2x_2=0\}$. Then we may similarly identify $\Gamma_\mm=\Gamma_1$ with $V/ \Lambda_K$, which is a torus of dimension $n-1=2$. The normalized Gr\"{o}ssencharacters ($\bmod\ 1$) form the dual of $\Gamma_1$, that is, a group isomorphic to $\Z^2$; all but one of  these are characters of infinite order. These characters take the angles in the complex embeddings of $K$ into account.

In order to give these Gr\"{o}ssencharacters  explicitly, we need a formula for the projection $\R^3\to V$ along the line $\R(1,1,0)\subset\R^3$, which corresponds to our fixed embedding of $\R^*_+$ into $\ak^*$. For this, take the dual basis $\{\bw_1,\bw_2\}$ to $\{\bv_1,\bv_2\}$ viewed in  the subspace $W$ orthogonal to $\R(1,1,0)$ with respect to the standard scalar product $\langle \bw, \bx\rangle= w_1 x_1 + w_2x_2+w_3x_3$:
\[
\bw_1 = (\frac{2}{3 \log \theta}, -\frac{2}{3\log \theta}, 0), \ \
\bw_2 = (-\frac{2\phi}{3 \log \theta}, \frac{ 2\phi}{3\log \theta}, \frac{1}{2 \pi}).
\]
Then the projection is given by $\bx\mapsto\langle \bw_1, \bx \rangle \bv_1+ \langle \bw_2, \bx \rangle\bv_2$. Now, given a nonzero ideal $\aaa=(\alpha)$ in $K$, by multiplying by $-1$ if necessary we may suppose that the generator $\alpha>0$. Then the algebraic conjugate of $\alpha$ in the  field $\Q(\vartheta^{+})$  takes the form  $\alpha^{+}=\sqrt{\frac{N(\aaa)}{\alpha}}e^{2\pi i\psi}$, where $\psi=\psi(\alpha)$ is well-defined $(\bmod \,1)$. The normalized Gr\"{o}ssencharacter $\gchar_{k_1, k_2}$ corresponding to $(k_1, k_2) \in \Z^2$ assigns to $\aaa$ the value
\[
 \gchar_{k_1, k_2}(\aaa) := e^{- 2 \pi i (k_1 \langle \bw_1, \bx \rangle + k_2 \langle \bw_2, \bx \rangle)},
\]
where  $\bx= \bx(\alpha) =( \log \alpha, - \frac{1}{2}\log \alpha+ \frac{1}{2}\log N(\aaa), 2 \pi \psi)$ is the image of $\alpha$ under the logarithm map.

\subsection{}\label{s4.3}
Theorem~\ref{tdistrib}, or rather Corollary~\ref{c1}, generalizes a result of Mitsui~\cite[Section~3]{M} on distribution of prime ideal numbers, which, in turn, generalizes earlier results of Hecke and Rademacher for quadratic number fields. In order to see, this consider the homomorphism $\theta\colon K^*_\infty\to \Kp_\infty$ that replaces every coordinate~$x_v$ of $x\in K^*_\infty$ for a real Archimedean place $v$ by $|x_v|$. Since the group $\Kp_\infty/\theta(U_{\mm,1})$ is divisible, the homomorphism
$\theta\colon K^*_\infty/U_{\mm,1}\to \Kp_\infty/\theta(U_{\mm,1})$ extends to a homomorphism $\theta_\mm\colon\ak^*/K^*W_\mm\to \Kp_\infty/\theta(U_{\mm,1})$. The homomorphism $\theta_\mm$ maps $\Gamma_{\mm}$ onto $\Kp_{\infty,1}/\theta(U_{\mm,1})$, where $\Kp_{\infty,1}=\Kp_\infty\cap K^*_{\infty,1}$. It follows that for $\rho_{\mm}: I_{\mm} \to \Gamma_{\mm}$,
 the image of the prime ideals under $\theta_\mm\circ\rho_\mm$, ordered by the norm, is equidistributed in $\Kp_{\infty,1}/\theta(U_{\mm,1})$. But we can conclude more.

Fix an ideal $\aaa\in I_\mm$ and consider its class $[\aaa]$ in the ray class group $C_\mm$. Identifying the set of cosets of~$K^*_{\infty,1}/U_{\mm,1}$ in $\Gamma_{\mm}$ with $C_\mm$, by Theorem~\ref{tdistrib} we conclude that the prime ideals in the class $[\aaa]$ have density $|C_\mm|^{-1}=h_K^{-1}\varphi(\mm)^{-1}$ in the set of all prime ideals, and their image under $\rho_\mm$, ordered by the norm, is equidistributed in the corresponding coset of $K^*_{\infty,1}/U_{\mm,1}$. The homomorphism~$\theta_\mm$ maps surjectively every such coset onto $\Kp_{\infty,1}/\theta(U_{\mm,1})$. It follows that: {\em  for each
class $[\aaa] \in C_{\mm}$,  the set $\{(\theta_\mm\circ\rho_\mm)(\pp)\mid\pp\in[\aaa]\}$, ordered by the norm, is equidistributed in $\Kp_{\infty,1}/\theta(U_{\mm,1})$}.

Consider now {\em ideal numbers} $\hat K^*$. They were introduced by Hecke, and are treated in Neukirch \cite[Section~VII.3]{Neu}. The construction can be described as follows. Put $\bar K_\infty=\prod_{v|\infty}\bar K_v$, where $\bar K_v$ is the algebraic closure of $K_v$, so $\bar K_\infty\cong\C^{r_1+r_2}$. Let~$\PP$ be the group of principal fractional ideals, and $I$ be the group of all fractional ideals.  Since the group~$\bar K^*_\infty$ is divisible, the canonical homomorphism $\PP=K^*/\OO^*\to \bar K^*_\infty/\OO^*$ can be (noncanonically) extended  to a homomorphism $I\to \bar K^*_\infty/\OO^*$. Consider the image $J\subset \bar K^*_\infty/\OO^*$ of~$I$ under this homomorphism, and then take $\hat K^*$ to be the preimage of $J$ in $\bar K^*_\infty$. We remark that Neukirch considers $\hat K^*$ as a subgroup of $K_\C=\prod_\tau\C\cong\C^n$, where $\tau$ runs over all embeddings $K\hookrightarrow\C$, rather than of $\bar K_\infty\cong \C^{r_1+r_2}$.

As before, consider an ideal $\aaa\in I_\mm$. Fix an ideal number $\hat a\in\hat K^*$ representing~$\aaa$. Then every ideal $\bb\in I_\mm$ in the class $[\aaa]\in C_\mm$ can be represented by a unique, up to a factor from~$U_{\mm,1}$, ideal number $\hat b$ lying in $K^*_{\mm,1}\hat a$; in other words, $\hat b\in\hat K^*$ is such that $\hat b\equiv\hat a\bmod\ \mm$. The homomorphism $\theta\colon K^*_\infty\to \Kp_\infty$ constructed earlier has an extension to~$\bar K^*_\infty$ defined in the same way. We therefore get an element $\theta(\hat b)\theta(U_{\mm,1})\in \Kp_{\infty}/\theta(U_{\mm,1})$ that depends on $\bb$ and~$\hat a$, but not on the choice of $\hat b$.  Finally, define $\rho(\bb):=s_\infty(N(\bb))^{-1}\theta(\hat b)\theta(U_{\mm,1})\in \Kp_{\infty,1}/\theta(U_{\mm,1})$. Then the result of Mitsui can be formulated by saying that for certain sets $A\subset \Kp_{\infty,1}/\theta(U_{\mm,1})$ with boundary of measure zero we have
$$
|\{\pp\in[\aaa]\mid N(\pp)\le x\ \text{and}\ \rho(\pp)\in A\}|=\frac{\tilde\lambda_\mm(A)x}{h_K\varphi(\mm)\log x}+o\left(\frac{x}{\log x}\right),
$$
where $\tilde\lambda_\mm$ is the normalized Haar measure on $\Kp_{\infty,1}/\theta(U_{\mm,1})$.
But by construction we have
$$
\rho(\pp)=\rho(\aaa)(\theta_\mm\circ\rho_\mm)(\aaa)(\theta_\mm\circ\rho_\mm)(\pp)^{-1},
$$
so this is the same equidistribution result as the one stated in italics above. Mitsui's result, however, is stronger in that he  obtains an estimate of the remainder term for sets $A$ of a particular shape, and obtaining
a remainder term is where the main work of Mitsui's paper  lies.

\subsection{}\label{s4.4}
As was already mentioned in the introduction, our ergodicity result is closely related to the analysis of Bost-Connes type systems.

A {\em $C^{*}$-dynamical system} consists of a  C$^{*}$-algebra of observables $A$ and a one-parameter family of automorphisms $(\sigma_t)_{t \in \R}$. A  {\em KMS-state}  at inverse temperature $\beta >0$ is a state $\varphi$ of $A$ which fulfills the KMS$_{\beta}$-condition with respect to $\sigma_t$, which says that for each $x,y \in A$ there exists a bounded holomorphic function $F_{x, y}(z)$ on the open strip $0 < \operatorname{Im}(z) < \beta$ which continuously extends to the boundary of the strip with boundary values satisfying
$$
F_{x,y}(t) =  \varphi (x \sigma_t(y))\ \ \text{and}\ \ F_{x,y}(t+i \beta) =  \varphi ( \sigma_t(y)x)
$$
for all $t\in\R$.

The Bost-Connes phase transition theorem \cite[Theorem 5]{bos-con} states that for each inverse temperature  $\beta>1$ the set of   KMS$_{\beta}$-states of the Bost-Connes C$^{*}$-dynamical system, attached to a certain Hecke pair, forms an infinite-dimensional  simplex whose extreme points are pure states $\varphi_{\beta, \chi}$ parametrized by complex embeddings $\chi\colon \Q^{ab} \to \C$, where $\Q^{ab}$ is the field generated by all roots of unity. These embeddings are described by injective characters $\chi\colon \Q/\Z \to \C$, and the restrictions of $\varphi_{\beta, \chi}$ to the ``diagonal'' subalgebra $C^{*}(\Q/\Z)\subset A$ are given by
$$
\varphi_{\beta, \chi}(\gamma) = \frac{1}{\zeta(\beta)} \sum_{n=1}^{\infty} \chi(\gamma)^n n^{-\beta},
$$
where $\zeta(\beta)$ denotes the Riemann zeta function. These states are type I$_{\infty}$ factor states. The system carries a $\Gal(\Q^{ab}/\Q)$-action and exhibits a symmetry breaking phase transition: for $0 < \beta \le 1$
the unique KMS$_{\beta}$-state is necessarily invariant under this action,
while for $\beta>1$ this action permutes the pure states. Finally, for $0 < \beta \le 1$  the unique  KMS$_{\beta}$-state is a type III$_1$ factor state.

To make the connection between our ergodicity result and this theorem more precise, let us  state our result in yet another form. Let
$$
C_K :=\ak^*/K^*
$$
 be the idele class group, and consider the balanced product $C_K\times_{\ohs}\akf$. The action $g(x,y)=(xg^{-1},gy)$ of $\akf^*$ on $\ak^*\times\akf$ defines an action of the group $I\cong\akf^*/\ohs$ of fractional ideals on $C_K\times_{\ohs}\akf$. Then a statement equivalent to Theorem~\ref{tmain1} is that the action of $I$ on $C_K\times_{\ohs}\akf$ is ergodic. Indeed, as in Section~\ref{s1}, ergodicity of this action is equivalent to ergodicity of the action of $I$ on $C_K\times_{\ohs}(\akf/W_\mm)=(C_K/W_\mm)\times_{\ohs}\akf$ for all $\mm$, which, in turn, is the same as ergodicity of the action of $I_\mm$ on $(C_K/W_\mm)\times(\aks/\oas^*)$, and this is exactly Theorem~\ref{tmain2}.

Assume now that $K$ is a number field and consider the compact group
$$
C^1_K :=\ker(N\colon C_K\to\R^*_+).
$$
Identifying it with the quotient $C_K/s(\R^*_+)$, where $s$ is as before a section of the norm map, we get an action of $I$ on $C^1_K\times_{\ohs}\akf$. Then ergodicity of the action of $I$ on $C_K\times_{\ohs}\akf$ amounts to the following: {\em the action of $I$ on $C^1_K\times_{\ohs}\akf$ is ergodic of type} III$_1$. For $K=\Q$ this  is known to be equivalent to the special case  of the Bost-Connes phase transition theorem \cite{bos-con} above for inverse temperature $\beta=1$. More precisely, ergodicity of the action of $\Q^*$ on $\A_{\Q}$ is equivalent to the statement that the unique $\Gal(\Q^{\rm{ab}}/\Q)$-invariant KMS$_1$-state on the Bost-Connes system is a type III$_1$ factor state. From this fact one easily deduces that this state is the only KMS$_1$-state on the Bost-Connes system, see the Remarks in~\cite{N1}.

For general number fields our ergodicity result is stronger than those proved  for Bost-Connes type systems~\cite{LLN,N2}, as  such systems have  the group $C^1_K$ replaced by its quotient group $\Gal(K^{ab}/K)\cong C_K/\overline{\Kp_\infty}$. The arguments of this paper can be adapted for all inverse temperatures $\beta\in(0,1]$. Here one proves ergodicity of the action with respect to a different measure $\mu_{ \beta}$  that depends on $\beta$, see~\cite{N1,LLN}. When complemented with the technique developed in~\cite{LLN} for $\beta>1$ this leads to the following generalization of the Bost-Connes phase transition theorem, formulated in dynamical systems terms.

\smallskip

\noindent{\bf Theorem.} {\it For a number field $K$ the following hold:
\begin{enumerate}
\item[(i)]
for every $\beta\in(0,1]$ there exists a unique Borel measure $\mu_\beta$ on $C^1_K\times_{\ohs}\akf$
with $\mu_\beta(C^1_K\times_{\ohs}\hat\OO)=1$ and with $\mu_\beta(\aaa\,\cdot)=N(\aaa)^{-\beta}\mu_\beta$ for all $\aaa\in I$; furthermore, the action of $I$ on $(C^1_K\times_{\ohs}\akf,\mu_\beta)$ is ergodic of type {\rm III}$_1$;
\item[(ii)]
for every $\beta\in(1,+\infty)$ there exists a bijective correspondence between points of $C^1_K\times_{\ohs}\ohs\cong C^1_K$ and extremal Borel measures $\nu$ on $C^1_K\times_{\ohs}\akf$ such that $\nu(C^1_K\times_{\ohs}\hat\OO)=1$ and $\nu(\aaa\,\cdot)=N(\aaa)^{-\beta}\nu$ for all $\aaa\in I$; namely, the measure corresponding to a point $x\in C^1_K\times_{\ohs}\ohs$ is given by $\zeta_K(\beta)^{-1}\sum_{\aaa\in I}N(\aaa)^{-\beta}\delta_{\aaa x}$.
\end{enumerate}
}

\smallskip

A similar result is true for function fields, but then instead of type III$_1$ in part (i) we get type III$_{q^{-\beta}}$, where~$q$~is the number of elements in the constant field of $K$.

\subsection{}\label{s5func}
In the global function field case with field of constants $\F_q$ the equidistribution result analogous to Theorem \ref{tdistrib} fails because the
$L$-function attached to the trivial Gr\"{o}ssencharacter has an infinite set of simple poles on the line
$\operatorname{Re}(s)=1$ at $s=1+\frac{2 \pi i k}{\log q}$.
These extra poles produce oscillations in the counting function of the Dirichlet series coefficients below $x$. For geometric extensions the pole locations and residues for  the trivial character do not change on the line $\operatorname{Re}(s)=1$, and the ergodicity result can still be  obtained in this case.
However for nongeometric extensions, the  change to the ground field $\F_{q^m}$ produces an $L$-function for the trivial character having  a closer spacing of poles on the line $\operatorname{Re}(s)=1$.

Assume $K$ is a global function field and $L$ is a finite abelian extension of $K$ unramified at all finite places $\pp\notin S(\mm)$. As in Theorem~\ref{tmain3}, consider the action of $I_\mm$ on $\Z\times\Gal(L/K)\times\aks/\oas^*$, given by
\[
\aaa(n,g,x)=(n+\deg\aaa,\sigma(\aaa)g,\aaa x).
\]
If the extension $L/K$ is not geometric, then this action is {\em not} ergodic, since the image $\Gamma$ of $I_\mm$ in $\Z\times\Gal(L/K)$ under
the homomorphism $\aaa\mapsto(\deg\aaa,\sigma(\aaa))$ is a proper subgroup. Namely, if the constant field of $L$ is $\F_{q^m}$, then $\Gamma$ is the kernel of the homomorphism $$\Z\times\Gal(L/K)\to\Gal(\F_{q^m}/\F_q),\ \ (n,g)\mapsto(g|_{\F_{q^m}})\operatorname{Frob}_{\F_{q^m}/\F_q}^{-n}.$$
Nevertheless Theorem~\ref{tmain3} can be formulated by saying that the kernel of the homomorphism $I_\mm\to\Z\times\Gal(L/K)$, $\aaa\mapsto(\deg\aaa,\sigma(\aaa))$, acts ergodically on $\aks/\oas^*$, and in this form the theorem remains true for nongeometric extensions. To put it differently, the asymptotic range of the cocycle $c$ defined in the proof of Theorem~\ref{tmain3} is exactly the subgroup $\Gamma\subset \Z\times\Gal(L/K)$. This is basically a reformulation of results of \cite[Section~2]{NR}, and this can also be proved using the same arguments as in our proof of Theorem~\ref{tmain3}.

\medskip

\paragraph{\bf Acknowledgments.} This work came about through
the authors' participation at  the  conference {\em Noncommutative Geometry: Multiple Connections}, held at
Ohio State University, May 2012. The conference was supported by the Mathematics Research Institute,
the NSF, and Ohio State University.

\bigskip

\end{document}